\documentclass{article}
\usepackage{amsmath}
\usepackage{latexsym}
\usepackage{amssymb}
\title{A selection principle\\in deformation
quantization}
\author{Murray Gerstenhaber\footnote{University of Pennsylvania,
Department of Mathematics, Philadelphia, PA 19104-6395;
\texttt{mgersten@math.upenn.edu}}}
\begin{document}
\maketitle
\newtheorem{theorem}{Theorem}
\newtheorem{corollary}{Corollary}
\renewcommand{\abstractname}{}

\begin{abstract}\noindent Deformation quantization
produces families of mathematically equivalent quantization
procedures from which one must select the physically meaningful
ones. As a selection principle we propose that the procedure must
allow enough `observable' energy distributions, i.e., ones for
which no pure quantum state will appear with negative probability
and must further have the property that for these the uncertainty
in the probability distribution of the quantum states must not
exceed that of the original distribution. For the simple harmonic
oscillator we show that this allows only the classic
Groenewold-Moyal (skew-symmetric) form.
\end{abstract}

The idea of negative probabilities is not new, going at least as
far back as Wigner and perhaps even into the 19th century. It has
been likened popularly to observers on the sidelines of a soccer
game seeing (non-negative) probability distributions of the
coordinates of the ball, none of which may seem unusual, but from
which they deduce that the probability distribution of the ball
over the interior of the entire playing field has points where it
is negative. Suppose, however, that we have a system for which we
know some distribution of energies. Deformation quantization
generally produces a family of `cohomologically equivalent'
quantizations of the system, each of which together with the
energy distribution assigns a probability, possibly negative, to
each pure quantum state. An energy distribution will be called
\emph{observable} with respect to a given quantization procedure
if each of these probabilities is in fact non-negative. We will
say that there are \emph{enough} observable distributions if every
distribution is in the closure of the linear space spanned by the
observable ones.

With any quantization procedure one can associate to the original
distribution of energies two measures of uncertainty (=standard
deviation), that of the original distribution and that of the
distribution of energies it produces in the various pure quantum
states. The latter will be called the \emph{quantum uncertainty};
The selection principle proposed here is that amongst
cohomologically equivalent quantizations only those are physically
meaningful for which there exist enough observable distributions,
and where for each observable distribution the quantum uncertainty
does not exceed that of the original distribution; briefly
\emph{quantization must decrease uncertainty}.

In the case of the simple harmonic oscillator, examining the
entire family of quantizations possible through deformation yields
identities involving Laguerre polynomials which, aside from the
present method of derivation, are generally not new. They show,
however, that for each of these quantizations there is a natural
infinite sequence of `basic' observable distributions from which
any pure state can be recovered as a linear combination. In
particular there are always enough observable distributions, so
for the simple harmonic oscillator this by itself is no
restriction on deformation quantization. However, the only
quantization that the inequality on uncertainties allows is the
Groenewold--Moyal form (although the normal or anti-normal form
may be meaningful when quantizing fields). To single out that form
in the case of the simple harmonic oscillator it would be
sufficient, as will be shown, to require that the quantum
uncertainty of an observable energy distribution approach that of
the original distribution as the energy tends to infinity; perhaps
this alone would be sufficient in general.

\section{Some basic algebraic deformation theory}
The seminal paper in deformation quantization is that of Bayen,
Flato, Fr{\o}nsdal, Lichnerowcz and Sternheimer  \cite{BFFLS}
1978, some essential ideas and results of which are simply
reproduced here without further attribution. This approach to
quantization has been exceedingly fruitful. Some of the subsequent
developments are summarized in \cite{BGGS:Quantum}, which is
complemented by an extensive bibliography. (A useful recent
introduction to the theory is the note of Hirshfeld and Henselder
\cite{HirshHens:Def}.) We begin with a very brief review of
algebraic deformation theory, introduced by the author in
\cite{G:DefI}.

Let $\mathcal{A}$ be an algebra which here (and generally in any
physical theory) will be assumed to be over the real or complex
numbers, e.g., the algebra of functions on a phase space, but in a
more general context could be over an arbitrary commutative,
unital ring. A \emph{deformation} of $\mathcal{A}$ is a new
associative `star' multiplication expressible as a formal power
series
\begin{equation} a*b= ab + \hbar C_1(a,b) + \hbar^2 C_2(a,b)
+\dots.\label{eq:basic} \end{equation} Here $\hbar$ is for the
moment just a formal parameter and the $C_i$ are bilinear maps
from $\mathcal{A} \times \mathcal{A}$ to $\mathcal{A}$. We should
like the star product to be defined on the same underlying vector
space as $\mathcal{A}$ but the introduction of the formal
parameter $\hbar$ generally makes it necessary to extend the
coefficients. Frequently the extension is to power series in
$\hbar$ and one views the functions $C_1, C_2, \dots$ as having
been extended to be bilinear not only over $\mathbb{R}$ or
$\mathbb{C}$ but also over these power series. Often one tacitly
assumes this done and in favorable cases the series actually
converge for sufficiently small values of the variable $\hbar$. In
deformation quantization, however, this extension of coefficients
is not appropriate because we are forced to consider power series
in $1/\hbar$. The correct extension of coefficients in this case
is to the field of Laurent series in $1/\hbar$, i.e.
$\mathbb{C}[[1/\hbar]][\hbar]$, but this puts additional
restrictions on the $C_i$ in order for the star multiplication to
be meaningful. It is sufficient that the star product
(\ref{eq:basic}) be \emph{locally finite}, i.e., that for any $a,
b \in \mathcal{A}$ only a finite number of the $C_i(a,b)$ be
non-zero. (This is likewise frequently unmentioned but it is
generally automatically satisfied when the $C_i$ are
bidifferential operators). The use of Laurent series in $1/\hbar$
rather than power series in $\hbar$ has, as we will see, a
profound effect on the structure of the resulting algebra.

If we set $C_0(a,b)=ab,$ the original associative multiplication,
then associativity of the star multiplication is equivalent to the
condition that
\begin{equation*}
\sum_{i+j=n}[C_i(C_j(a,b),c)-C_i(a,C_j(b,c))] = 0, \quad\text{all
}\, n\ge 1;\, i,j \ge 0\end{equation*} for all $a,b,c$ in the
algebra $\mathcal{A}$. Transposing to the right side all terms
with $i$ or $j$ equal to 0, this becomes
\begin{equation}\label{eq:assoc2}\begin{aligned}
\sum_{i+j=n;\, i,j >0}[&C_i(C_j(a,b),c)-C_i(a,C_j(b,c))] =\\&
aC_n(b,c)-C_n(ab,c)+ C_n(a,bc)-C_n(a,b)c, \quad\text{all }\, n\ge
1.\end{aligned}\end{equation} These are generally difficult
conditions to meet. In the Hochschild cohomology theory, each
$C_i$ is a 2-cochain of $\mathcal{A}$ with coefficients in itself,
and the right side in (\ref{eq:assoc2}) is the Hochschild
coboundary of $C_n$. For $n = 1$ the left side is zero, so the
coboundary of $C_1$ is zero, that is, $C_1$ is a 2-cocycle. This
(or more properly, its cohomology class) is often called the
\emph{infinitesimal} of the deformation. A basic problem, given an
infinitesimal deformation, is to construct a deformation which has
it for infinitesimal. For $n=2$ the left side of (\ref{eq:assoc2})
is something constructed from $C_1$ which in fact will always be a
3-cocycle in the Hochschild theory, and the first requirement is
that it be a coboundary, namely the coboundary of $C_2$. Having
$C_1$ and $C_2$, the left side of (\ref{eq:assoc2}) with $n=3$
will again be a cocycle which is required to be a coboundary, and
so forth for all $n$. Unless we have some control of the
Hochschild cohomology of $\mathcal{A}$ it is clear that the
construction of deformations in such a step-by-step manner will
not be easy.

Fortunately, there is one case in which all the conditions for
associativity are automatically satisfied. A \emph{derivation} $D$
of $\mathcal{A}$ is a linear mapping of $\mathcal{A}$ into itself
such that $D(ab) = (Da)b + a(Db)$. We have been careful with the
order of the variables $a$ and $b$ here because the multiplication
in the algebra $\mathcal{A}$ need not have been commutative,
although that is the case for any usual algebra of functions.
Ordinary differentiation in the algebra of infinitely
differentiable functions on $\mathbb{R}$ is a derivation. Suppose
that $D'$ and $D''$ are commuting derivations of $\mathcal{A}$.
Then the star multiplication defined by
$$a*b = ab +\hbar D'a\,D''b + \frac{\hbar^2}{2!}D'^2a\,D''^2b +
\frac{\hbar^3}{3!}D'^3a\,D''^3b + \dots$$ will be associative, cf.
\cite{G:DefI}, something easily verified by direct computation.
This (despite the consternation of mathematicians) is frequently
written as
\begin{equation} a*b =
a\,e^{\hbar\stackrel{\leftarrow}{D'}\stackrel{\rightarrow}{D''}}b.
\label{eq:normal} \end{equation} In the star product of
(\ref{eq:basic}) the $C_n$ are then given by $C_n(a,b) =
(1/n!)D'^na\,D''^nb.$ More generally, if $D'_1,\dots, D'_r,
D''_1,\dots, D''_r$ are all mutually commuting derivations then
\begin{equation*}
a*b = a\,e^{\hbar\sum_{i=1}^r
\stackrel{\leftarrow}{D'_i}\stackrel{\rightarrow}{D''_i}}b
\label{eq:normal2}
\end{equation*}
is again an associative multiplication.

A basic example illustrating (\ref{eq:normal}) is that where
$\mathcal{A} = \mathbb{C}[q,p], D'= \partial_q, D'' =
i\partial_p.$ Then $q*p = qp +ih$, while  $p*q =pq,$ so $[q,p]_* =
q*p-p*q = i\hbar.$ \emph{The last equation determines the
structure of the deformed algebra up to isomorphism but not the
deformation} (which contains more information), as different
deformations can give isomorphic algebras. In fact, referring to
the basic equation (\ref{eq:basic}), suppose that $T:\mathcal{A}
\to \mathcal{A}$ is a linear map of the underlying vector space
onto itself of the form $Ta = a+ \hbar \tau_1a +\hbar^2\tau_2a +
\dots.$ Defining $a*'b = T^{-1}(Ta\,*\,Tb)$ and denoting by
$\mathcal{A}_*, \mathcal{A}_{*'}$ the algebras with these two
multiplications, the map $T$ is an isomorphism $\mathcal{A}_{*'}
\to \mathcal{A}_*$, so the new multiplication is also associative.
We say that the deformations given by $*$ and $*'$ are
\emph{cohomologically equivalent} (c-equivalent), the adjective
emphasizing that despite the isomorphism there may be some
physical differences between the results. In particular, consider
(\ref{eq:normal}) and take $T = T_{\lambda} =
e^{\lambda\hbar\,D'D''}$ where $D'D''$ is just the composite of
the two derivations and $\lambda$ is an arbitrary constant. Then
it is easy to check that the resulting deformation is given by
\begin{equation} a*_{\lambda}b =
a\,e^{\hbar\stackrel{\leftarrow}{D'}\stackrel{\rightarrow}{D''} -
\lambda\hbar(\stackrel{\leftarrow}{D'}\stackrel{\rightarrow}{D''}
- \stackrel{\leftarrow}{D''}\stackrel{\rightarrow}{D'})}\,b
\label{eq:lambdaMoyal} \end{equation} In particular, for $\lambda
= 1/2$ we have
\begin{equation} a\,*_{1/2}\,b =
a\,e^{\frac{\hbar}{2}(\stackrel{\leftarrow}{D'}\stackrel{\rightarrow}{D''}-
\stackrel{\leftarrow}{D''}\stackrel{\rightarrow}{D'})}\,b.
\label{eq:skewform} \end{equation} In the context of quantum
theory, the deformation given by the commuting derivations
$D',D''$ in (\ref{eq:normal}) is generally called the `normal'
form, that obtained by interchanging $D'$ and $D''$ or by setting
$\lambda = 1$ in (\ref{eq:lambdaMoyal}) is the `anti-normal' form,
and the skew symmetric form of (\ref{eq:skewform}) is the
`Groenewold-Moyal' form (GM)\footnote{While often attributed
solely to Moyal, the basic idea is present earlier in the work of
Groenewold, and anticipated even earlier in works of Wigner and
Weyl.} \cite{Groen:Principles, Moyal}. Historically, the fact that
cohomologically equivalent deformations may not be physically
equivalent immediately raised the problem of selecting from a
family of c-equivalent deformations those which are physically
meaningful. This was already addressed in \cite{BFFLS} where one
important reason cited for preferring the GM form is its greater
symmetries. Another is homological: Every 2-cocycle of
$\mathbb{C}[q,p]$ with coefficients in itself can be written
uniquely as a sum of a symmetric part and a skew part. Both parts
are again cocycles, but the symmetric part is always a coboundary;
there is, up to constant multiples, a unique skew 2-cocycle, and
that cocycle is a biderivation, i.e., a derivation as a function
of each argument. As mentioned, in the case of the simple harmonic
oscillator our selection principle allows only the GM form, but
like the foregoing principles (symmetry, cohomological uniqueness)
our selection principle should apply to many other cases. (The
normal and anti-normal forms are excluded in the case of a single
oscillator but may be essential when one has infinitely many, as
when quantizing a field, and one must normalize the lowest energy
level of each to zero to avoid having that of the whole be
infinite.)

Consider now the choice of coefficients of a deformed algebra. If
we do not have some specific information about the 2-cochains
$C_i$ in (\ref{eq:basic}) other than that they give a deformation
(or about the derivations $D'$ and $D''$ in (\ref{eq:normal})
other than that they commute), then for coefficients one must take
the power series ring $\mathbb{C}[[\hbar]]$, else the formulas
will not be meaningful. In this classic approach some basic
algebraic properties of $\mathcal{A}$ are preserved. In
particular, if a non-zero element $a$ of $\mathcal{A}$ is not a
zero divisor in the original multiplication (i.e., if there is no
$b \ne 0$ such that either $ab = 0$ or $ba = 0$) then $a$ will not
be a zero divisor in the deformed algebra. Similarly, if $a$ was
invertible then it will continue to be so. It follows that a
deformation of an integral domain will continue to be an integral
domain, and a deformation of a division ring (= skew field) will
again be a division ring. A deformation of a unital algebra
remains unital and the deformation will be, in fact, c-equivalent
to one in which the original unit remains the unit. However, a
deformation of a commutative algebra like $\mathbb{C}[q,p]$ need
not remain commutative; this is the basis of quantization.

As mentioned, in a classical algebraic deformation we generally
hope that the power series which are encountered actually converge
for sufficiently small values of the deformation parameter
$\hbar$, but for purely algebraic purposes this may not be
necessary. Suppose now, however, that the deformation has the
local finiteness property that for every $a$ and $b$ in the
original (undeformed) algebra there is an $N$ such that $C_i(a,b)
= 0$ for all $i > N$. The deformed algebra will then already be
defined over the polynomial ring $\mathbb{C}[\hbar]$, and we can
extend coefficients, if we wish, to the field of Laurent series in
$1/\hbar$. This is the case, for example, with $\mathcal{A} =
\mathbb{C}[q,p]$ and $D' = \partial_q =
\partial/\partial q, D'' = i\partial_p$ in (\ref{eq:normal}). We
must now also be careful in the definition of cohomological
equivalence to require that $T$ also be locally finite, i.e., that
for all $a$ there is an $N$ such that $\tau_i(a) = 0$ for $i > N$.
That is certainly the case for the $T$ which gives the equivalence
between the normal and Groenewold--Moyal deformations of
$\mathbb{C}[q,p]$ (with $D' = \partial_q, D'' = \partial_p$).

When, in the locally finite case, we extend coefficients to
Laurent series in $1/\hbar$ the structure of the deformed algebra
may be very different from that obtained with power series as
coefficients. What was before deformation an integral domain may
acquire infinitely many orthogonal idempotents; as a result there
may be no natural way to apply a contraction in the sense of
{{\.{I}}n{\"{o}}n{\"{u}} and Wigner \cite{InWig:Contraction} to
recover the original algebra. While this is inherent in
deformation quantization, it raises difficult purely algebraic
questions about the structure of the deformed algebra. The same is
true, of course, of all algebras obtained by c-equivalent
deformations since they are all algebraically isomorphic.

\section{Deformation quantization of the simple harmonic
oscillator} The foundational paper \cite{BFFLS} showed, in
particular, that quantization of the simple harmonic oscillator
could be viewed as an exercise in the deformation of the
polynomial ring $\mathbb{C}[q,p]$, where now  $q$ and $p$ are
viewed as the position and momentum coordinates on the phase space
$\mathbb{R}^2$. The Hamiltonian function for the simple harmonic
oscillator is
\begin{equation}H(q,p) = \frac{p^2}{2m} +
\frac{m\omega^2}{2}q^2.\label{eq:Ham}\end{equation}
 In the classic
approach to quantization one substitutes for $q$ and $p$ operators
$Q$ and $P$ which satisfy the fundamental commutation relation
$[Q,P] = i\hbar$. Generally this involves some ambiguity, for the
Hamiltonian may contain monomials of positive degree
simultaneously in $p$ and $q$, but that is not a problem here. We
can take, e.g., $Q=$multiplication by $q$, $P = -i\hbar
\partial_q$ and with this seek solutions to the Schr\"odinger
equation
\begin{equation*}\label{eq:Sch} i\hbar \dot{\psi} = H\,\psi.\end{equation*}
This in effect chooses a specific \emph{representation} of the
Weyl algebra $\mathbb{C}\{q,p\}/(qp-pq-i\hbar)$ and with this
choice $\psi$ is viewed as a function of $q$ and $t$. (The present
$H$ is time independent.) Mathematically this does not yet
introduce any quantization; the latter is forced by the physical
requirement that $\psi$ be square integrable with absolute value
tend to zero at $\pm\infty.$ By contrast, \emph{the deformation
approach chooses a deformation which gives rise to the Weyl
algebra}, e.g. that in (\ref{eq:normal}) (i.e., such that the
commutator of $q$ and $p$ is essentially their Poisson bracket)
and rewrites the Schr\"odinger equation in the form
\begin{equation} i\hbar \dot{\psi} = H*\psi.
\label{eq:defSch}\end{equation} This will be called the
``deformation-Schr\"odinger" or d-Schr\"odinger equation.

Were $H$ a matrix operating on a vector $\psi$, the solution would
be $e^{-iHt/\hbar}\psi(0)$. Here, bearing in mind that $H$ is now
an element of a non-commutative algebra with multiplication $*$
one must still compute the exponential
$\mathrm{exp}_*(-iHt/\hbar)$ where $\mathrm{exp}_*$ indicates that
the exponential must be computed using the deformed
multiplication. The problem is to express the result, which is an
element of the underlying vector space of the original
(undeformed) algebra of functions on phase space, without
reference to the deformed multiplication. There are now different
possible choices for the deformed multiplication $*$ but
$\mathrm{exp}_*(-iHt/\hbar)\psi(0)$ will always be a solution to
the d-Schr\"odinger equation.

In deformation quantization generally, one knows from \cite{BFFLS}
that
\begin{equation}\mathrm{exp}_*(-itH/\hbar)=
\sum_E\pi_E\,e^{-itE/\hbar} \label{eq:energy}\end{equation} where
the sum in the Fourier--Dirichlet series on the right is over the
allowable energy levels $E$ and the $\pi_E$ are functions on the
phase space which are orthogonal idempotents in the *
multiplication whose sum is 1. One has $H*\pi_E = E\,\pi_E.$ (Note
that we have tacitly extended coefficients to Laurent series in
$1/\hbar$ and this has introduced zero-divisors into the new
algebra.)  Further, these idempotent functions when integrated
over all of phase space will yield a common constant which, in the
case of simple harmonic motion whose phase space is the $q,p$
plane, is $2\pi\hbar.$ We should like to interpret the $\pi_E$,
which are functions on phase space, as giving a probability
distribution there but in general they may take on negative values
(depending on the deformation chosen). Although negative
probabilities may not be directly observable, we shall see that at
least in the case of simple harmonic motion a reasonable
interpretation as probabilities may still be possible.

Following an idea often credited to Dirac, it is convenient to
transform the Hamiltonian (\ref{eq:Ham}) into ``holomorphic
coordinates'' by setting
\begin{equation*} a =\sqrt{\frac{m\omega}{2}}(q+i\frac{p}{mw}), \quad
\bar{a} =
\sqrt{\frac{m\omega}{2}}(q-i\frac{p}{mw}).\label{eq:holo}\end{equation*}
With this one has $$H=\omega a\bar{a}.$$

The simplest deformation of $\mathbb{C}[a,\bar{a}]$ one can now
choose is the normal form defined by setting
\begin{equation*}
f\,*_N\,g=f\,
e^{\hbar\stackrel{\leftarrow}{\partial_a}\stackrel{\rightarrow}{\partial_{\bar{a}}}}
\,g.
\end{equation*}
One then has $[a,\bar{a}]_{*_N} = \hbar$, which is equivalent to
$[q,p]_{*_N} = i\hbar.$  With this quantization we must compute
$\mathrm{exp}_{*N}(-iHt/\hbar)= \mathrm{exp}_{*_N}(-i\omega
ta\bar{a}/\hbar)$. The d-Schr\"odinger equation (\ref{eq:defSch})
actually is a simple first order partial differential equation
which shows, in particular, that the solution is a function of
$a\bar{a}$ only. Writing $a\bar{a}=s$ and denoting the solution by
$F(s,t)$, the d-Schr\"odinger equation becomes
$$i\hbar \partial_tF(s,t) = \omega sF(s,t) +
\omega\hbar\partial_sF(s,t).$$ The required solution, which must
have the value 1 at $t=0$, is $F(s,t) =
e^{-s/\hbar}\mathrm{exp}(e^{-i\omega t}s/\hbar)$, so we have
\begin{equation*}\mathrm{exp}_{*_N}(-i\omega
ta\bar{a}/\hbar)= e^{-a\bar{a}/\hbar}\mathrm{exp}(e^{-i\omega
t}a\bar{a}/\hbar), \end{equation*} where on the right one has the
\emph{ordinary} exponential. Expanding the expression on the right
and writing $a\bar{a}/\hbar = H/\hbar\omega = \mu,$ the
coefficient of $e^{-in\omega t}$ becomes $e^{-\mu}{\mu}^n/n!$.
Comparing with (\ref{eq:energy}) we see that the allowable values
for the energy are $E_n = n\hbar\omega$ and the corresponding
$\pi_n^{(N)} = e^{-\mu}{\mu}^n/n!$ (where $N$ indicates that the
normal form is used). We will write the $\pi_n^{(N)}$ (and
generally those which arise with any quantization) as functions of
$H$ with the latter viewed as the energy function on phase space.
For any value of $H$ the $\pi_n^{(N)}$ are just the terms in the
classical Poisson distribution with mean $\mu = H/\hbar\omega.$
This probability distribution is sometimes called the `law of rare
events'\footnote{Curiously, it arose neither from gambling nor
physics but from Sim\'eon-Denis Poisson's study of the French
judicial process in his ``Recherches sur la probabilit\'e des
jugements en mati\`ere criminelle et mati\`ere civile", 1837} : if
$\mu$ is the mean number of events seen in unit time (or space)
then the probability that in a given unit of time (or space) one
will see exactly $n$ events is $\pi_n^{(N)} = e^{-\mu}{\mu}^n/n!.$
For example, if misprints are rare and the average number on a
page is $\mu$ then the probability that a page will contain
exactly $n$ misprints is $e^{-\mu}{\mu}^n/n!$ (but $\mu$ may vary
with the author). Unlike a Gaussian distribution, which depends on
two parameters, its mean and standard deviation (`uncertainty' in
physical terms, square root of its variance), the Poisson has but
one, its mean. The variance is identical with its mean and the
standard deviation is the square root of its mean. For large
values of the mean, the distribution resembles a Gaussian with
mean $\mu$ and standard deviation $\sqrt{\mu}$.

The question is how to interpret the appearance of the Poisson
distribution here, bearing in mind that we have somehow the
`wrong' quantization (or at least not that in textbooks, since the
lowest allowable energy is precisely zero, not $\hbar\omega/2$).
We will see that this quantization essentially presumes that we
know the mean energy of the oscillator precisely, something which
is not physically possible. In that impossible case it seems to
say that if the mean energy is $\mu$ then the probability that the
oscillator will be observed in a state with the quantum number $n$
is $\pi_n^{(N)} = e^{-\mu}{\mu}^n/n!.$ But note that even though
we have the `wrong' minimal energy, the spectrum here is simply
shifted by $\hbar\omega/2$ from the textbook case, so the
differences between allowable energy levels, which determine the
spectrum, coincide with the usual. Nevertheless, the normal form
of quantization is excluded by our selection principle. For with
it all true probability distributions of energy are observable
including a delta function, which has zero uncertainty, while the
quantum uncertainty is always positive. (In fact, we will see that
when the normal form is viewed as a limit the initial energy
distribution is a delta function.)

It is easy to verify that the $\pi_n^{(N)}$ are orthogonal
idempotents in the $*_N$ multiplication, summing to 1 and having a
fixed common integral over phase space: Direct computation shows
that the integral of each is $2\pi\hbar$ independent of $n$. That
$\sum_{n=0}^{\infty}\pi_n^{(N)} = 1$ follows simply from setting
$t=0$. Finally, one way to see that the $\pi_n^{(N)}$ are
orthogonal idempotents is to observe that although $*_N$ is a
non-commutative multiplication, the$\pi_n^{(N)}$ and
$\mathrm{exp}_{*_N}(-i\omega ta\bar{a}/\hbar)$ are all functions
only of the single element $a\bar{a}$ and hence all commute.
Comparing the expansions of the two sides of the equation
$\mathrm{exp}_{*_N}(-i\omega ta\bar{a}/\hbar)\,*_N\,
\mathrm{exp}_{*_N}(-i\omega ta\bar{a}/\hbar) =
\mathrm{exp}_{*_N}(-2i\omega ta\bar{a}/\hbar)$ will show that
$\pi_0^{(N)}$ is idempotent. Denoting it for the moment by $e$, we
clearly have for any idempotent $e$ that $e(1-e)=0$, and that
$1-e$ is again idempotent. Proceeding by induction will show that
the $\pi_n^{(N)}$ are all mutually orthogonal idempotents.

While the normal form of deformation has been excluded by our
selection principle, it already raises an interesting algebraic
question equally meaningful for all c-equivalent deformations. For
the moment, let $\mathcal{A}$ denote the polynomial ring
$\mathbb{C}[q,p]$ and let $\mathcal{A}_{\hbar}$ denote the algebra
to which we have deformed it. Note that as long as coefficients
are restricted to polynomials in $\hbar$ it is meaningful to let
$\hbar \to 0$ in order to recover the original algebra; this gives
the (only) correct statement of the correspondence principle.  The
resulting algebra is essentially the first Weyl algebra, a simple
algebra (i.e., one without proper two-sided ideals) whose
cohomology with coefficients in itself vanishes in all positive
dimensions, in particular in dimension 2, and which is therefore
rigid. With coefficients extended to
$\mathbb{C}[[1/\hbar]][\hbar]$ recovery of the original algebra by
letting $\hbar \to 0$ is no longer possible. We have an algebra in
which the identity has decomposed into a direct sum of infinitely
many orthogonal idempotents $e_1, e_2, \dots$ and whose precise
structure we no longer know. The theorem that there is no
degeneracy in one dimension suggests that each
$e_i\mathcal{A}_he_i$ has dimension 1 and that the same is
probably true for all $e_i\mathcal{A}_he_j$. The simplest
conjecture concerning structure would be that in each
$e_i\mathcal{A}_he_j$ we can choose an element $e_{ij}$ with
$e_{ii}$ = $e_i$ and $e_{ij}e_{jk} = e_{ik}$, and that the algebra
consists of linear combinations of these (but what beside the
finite ones may be allowed is not clear). Again, the same question
arises for all c-equivalent deformations of $\mathbb{C}[q,p]$.

\section{Quantizations with $\lambda \ne 0$} Following the
prescription in \S\,1, we now set $T_{\lambda} =
e^{\lambda\hbar\,D'D''}$ and define
\begin{equation*} f\,*_{\lambda}\,g =
T_{\lambda}^{-1}(T_{\lambda}f\,*_N\,T_{\lambda}\,g) =
f\,e^{\hbar((1-\lambda)\stackrel{\leftarrow}{\partial_a}\stackrel{\rightarrow}
{\partial_{\bar{a}}} - \lambda\stackrel{\leftarrow}
{\partial_{\bar{a}}}\stackrel{\rightarrow}{\partial_a})}\,g .
\end{equation*} There are now several approaches to computing
$\mathrm{exp}_{*_{\lambda}}(-iHt/\hbar)$. One can use the first
equality above to get expressions for the quantities which we will
now denote by $\pi_n^{(\lambda)}$, where $\pi_n^{(0)} =
\pi_n^{(N)}$; this will give them all in the form of power series.
The second approach, which we adopt, is to use only the value
obtained for $\pi_0^{(\lambda)}$ from the first method and then to
adapt the procedure in \cite[Appendix]{HirshHens:Def} to compute
$\mathrm{exp}_{*_{\lambda}}(-iHt/\hbar).$  The third approach is
to note, as in the preceding section, that
$\mathrm{exp}_{*_{\lambda}}(-iHt/\hbar)$ will be a function only
of $a\bar{a}$; denoting this again by $s$ and the result by
$F(s,t)$ one can solve the partial differential equation (now of
second order in $s$) which $F$ satisfies.  However, we will see
that one can also effectively solve the differential equation in
closed form once closed expressions for the $\pi_n^{(\lambda)}$
have been obtained from the second method by using the generating
function for the Laguerre polynomials. Since $T_{\lambda}a\bar{a}
= a\bar{a} + \hbar\lambda,$ with the first approach we have
\begin{equation*}\begin{split}
\mathrm{exp}_{*_{\lambda}}(-iHt/\hbar) &=
T_{\lambda}^{-1}\mathrm{exp}_{*_N}(-i\omega
tT_{\lambda}(a\bar{a}/\hbar))\\  &=
e^{-\lambda\hbar\partial_a\partial_{\bar{a}}}\mathrm{exp}_{*_N}
(-i\omega t(a\bar{a}/\hbar + \lambda)) \\
 &= e^{-\lambda i\omega t/\hbar}
e^{-\lambda\partial_a\partial_{\bar{a}}}e^{-a\bar{a}/\hbar}
\mathrm{exp}(e^{-i\omega t}a\bar{a}/\hbar) \\
&= e^{-\lambda i\omega t/\hbar}
e^{-\lambda\partial_a\partial_{\bar{a}}} \mathrm{exp}((e^{-i\omega
t}-1)a\bar{a}/\hbar).\end{split}\end{equation*} A simple
computation using the definition of the Laguerre polynomials then
shows that
\begin{equation*}
\mathrm{exp}_{*_{\lambda}}(-iHt/\hbar) = e^{-\lambda i\omega
t}\sum_{k=0}^{\infty}(-\lambda)^k(e^{-i\omega
t}-1)^k\mathrm{L}_k(a\bar{a}/\lambda\hbar)
\end{equation*} where $\mathrm{L}_k$ is the $k$th Laguerre
polynomial. We may now replace $a\bar{a}/\hbar$ on the right with
$H/\hbar\omega$ since the computations on the right now all take
place in the undeformed algebra. Expanding the right side then
gives
\begin{equation*}
\mathrm{exp}_{*_{\lambda}}(-iHt/\hbar) =
\sum_{n=0}^{\infty}(-1)^n\sum_{k=0}^{\infty}\lambda^k\binom{k}{n}
\mathrm{L}_k(H/\lambda\hbar\omega)e^{-i(n+\lambda)\omega t}
.\end{equation*} The spectrum has thus been shifted, the allowable
values of the energy are now $E=(n+\lambda)\hbar\omega$, and we
have
\begin{equation*}
\pi_n^{(\lambda)} =
(-1)^n\sum_{k=0}^{\infty}\lambda^{n+k}\binom{n+k}{k}\mathrm{L}_{n+k}(H/\lambda\hbar\omega),
\end{equation*} which can be negative for some $H$.
Using the generating function for the Laguerre polynomials,
\begin{equation*}
\frac{1}{1+x}\mathrm{exp}(\frac{zx}{1+x})
=\sum_{k=0}^{\infty}x^k(-1)^k\mathrm{L}_k(z),
\end{equation*}
the special case $n=0$ gives
\begin{equation*}
\pi_0^{(\lambda)} =
\sum_{k=0}^{\infty}\lambda^k\,\mathrm{L}_k\Bigl(\frac{H}{\lambda\hbar\omega}\Bigr)
=
\frac{1}{1-\lambda}\mathrm{exp}\Bigl(-\frac{H}{(1-\lambda)\hbar\omega}\Bigr).
\end{equation*} We can now adapt the procedure in
  \cite[Appendix]{HirshHens:Def} to get from this the  closed
form for all the $\pi_n^{(\lambda)}$, namely
\begin{equation*}
\pi_n^{(\lambda)} =
\frac{1}{1-\lambda}\Bigl(\frac{-\lambda}{1-\lambda}\Bigr)^n\,
\mathrm{L}_n\Bigl(\frac{H}{\lambda(1-\lambda)\hbar\omega}\Bigr)
\mathrm{exp}\Bigl(-\frac{H}{(1-\lambda)\hbar\omega}\Bigr).
\end{equation*} Letting $\lambda \to 0$ recovers original Poisson distribution.
Comparing the two expressions for $\pi_n^{(\lambda)}$ gives the
following identity involving Laguerre polynomials:
\begin{equation}\label{eq:LagId}\begin{split}
(-1)^n\sum_{k=0}^{\infty}\lambda^{n+k}\binom{n+k}{k}\mathrm{L}_{n+k}(\frac{z}{\lambda})
\qquad\qquad\qquad\qquad\hspace{.75in}\\
\qquad\qquad=\frac{1}{1-\lambda}\Bigl(\frac{-\lambda}{1-\lambda}\Bigr)^n\,
\mathrm{L}_n\Bigl(\frac{z}{\lambda(1-\lambda)}\Bigr)
\mathrm{exp}\Bigl(-\frac{z}{(1-\lambda)}\Bigr).
\end{split}
\end{equation}
(Setting $\lambda = 1/2$ and $z=0$ this asserts, for example, that
$\sum_k(1/2)^{n+k}\binom{n+k}{k} = 2$ independent of $n$,
something easily verified directly.) Multiplying the left side by
$e^{-z/\lambda}\mathrm{L}_m(z/\lambda) = e^{-z/\lambda}
\mathrm{L}_m((1-\lambda)z/\lambda(1-\lambda))$ and integrating,
the orthogonality relations of the Laguerre polynomials together
with an obvious change of variables gives
\begin{equation*}\int_0^{\infty}\mathrm{L}_m((1-\lambda)z)\mathrm{L}_n(z)e^{-z}\,dz
= \left\{\begin{array}{cl}
\binom{m}{n}\lambda^m\bigl(\frac{1-\lambda}{\lambda}\bigr)^n &
\mbox{ if $m\ge n$}\\0 & \mbox{otherwise} \end{array}. \right.
\end{equation*}
Replacing $\lambda$ by $1-\lambda$ and comparing the coefficients
of the powers of $\lambda$ on the two sides gives the coefficients
in the Fourier-Laguerre expansion of $z^k$:
\begin{equation}\int_0^{\infty}z^k\mathrm{L}_n(z)e^{-z}\,dz = \left\{
\begin{array}{cl}(-1)^n\binom{k}{n}k! &
\mbox{if $k \ge n$}\\0 & \mbox{ if $k < n$}
\end{array}.\right.\label{eq:fund}\end{equation}
This fundamental result, which will show the existence of enough
observable distributions, can also be derived in an elementary
way, since it is just the formula for the change of basis from the
Laguerre polynomials to the powers of $x$ in the the vector space
of polynomials in $x$: Let $A$ be the infinite lower triangular
matrix with rows and columns indexed by $0,1,2,\dots$ and $(i,j)$
entry equal to $(-1)^j\binom{i}{j}$, and $D$ be the infinite
diagonal matrix with diagonal entries $1/n!,\, n = 0,  1, 2,
\dots\,$\,. Letting $X$ be the infinite column vector
$(1,x,x^2,x^3, \dots)^t$, the $n$th entry in the vector $L = ADX$
is just the Laguerre polynomial $\mathrm{L}_n(x)$. To write $X$ in
terms of $L$ it is sufficient therefore to invert $AD$, the only
problem being the inversion of $A$. However, $A$ is equal to its
own inverse, for writing $x^r = (1-(1-x))^r = \sum
(-1)^i\binom{r}{i}(1-x)^i$ gives
$$\sum_i(-1)^i\binom{r}{i}(-1)^j\binom{i}{j} =
\left\{\begin{array} {cl} 1 & \mbox{if $j=r$}\\0& \mbox{otherwise}
\end{array}.\right.$$ Therefore $X = D^{-1}AL,$ which is
precisely what  (\ref{eq:fund}) asserts.

While $k$ is an integer in (\ref{eq:fund}), one can deduce more
generally that
\begin{equation*} \int_0^{\infty}z^p\mathrm{L}_n(z)e^{-z} =
(-1)^n\frac{\Gamma(p+1)^2}{n!\Gamma(p-n+1)}.\end{equation*} With
(\ref{eq:LagId}), the generating function for the Laguerre
polynomials gives the following closed form for the $*_{\lambda}$
exponential of $-iHt/\hbar$,
\begin{equation*}
\mathrm{exp}_{*_{\lambda}}(-iHt/\hbar)\quad  =
\quad\frac{e^{-i\lambda\omega t}}{1-\lambda + \lambda e^{-i\omega
t}}\,\mathrm{exp}\Bigr(\frac{2H}{\hbar\omega}\cdot
\frac{e^{-i\omega t}-1}{1-\lambda +\lambda e^{-i\omega t}}\Bigl);
\end{equation*}
at $\lambda = 1/2$ one recovers the known formula for the
Groenewold-Moyal case,
\begin{equation*}
\mathrm{exp}_{*_{1/2}}(-iHt/\hbar)\quad  =
\quad
\frac{1}{\cos(\omega
t/2)}\,\mathrm{exp}\Bigr(\frac{2H}{-i\hbar\omega}\,\tan(\omega
t/2)\Bigl).
\end{equation*}
The only values of $\lambda$ that need to be considered are
$0\le\lambda\le 1/2$. Setting $t=0$ shows that
$\sum\,\pi_n^{(\lambda)}= 1$ for any value of $\lambda$; that the
$\pi_n^{(\lambda)}$ are mutually orthogonal idempotents follows
exactly as in the case $\lambda = 0$.  We could, of course also
deduce this from the fact that the $\pi_n^{(\lambda)}$ are the
transforms of the $\pi_n^N$ by $T_{\lambda}$, but writing, as
before, $\mu = H/\hbar\omega$ we also have the duality relation
$$\int_0^{\infty}\pi_n^{(\lambda)}(\mu)\,\pi_m^{(1-\lambda)}(\mu)\,d\mu =
\delta_{n,m}.$$ The Groenewold-Moyal case ($\lambda = 1/2$) is
 self-dual.

\section{Negative probabilities and basic observable distributions}
Returning to the question of negative probabilities, unlike the
$\pi_n^{(N)}$, we can not view the $\pi_n^{(\lambda)}$ as giving
an ordinary probability distribution over the energy values $H$
since for all $n
> 0$ and any positive $\lambda$, $\pi_n^{(\lambda)}$  will be
negative for some positive value of $H$. (Moreover, we conjecture
that for any fixed $\lambda \ne 0,1$ and arbitrary  $H
> 0$ there must be an $n
> 0$ such that $\pi_n^{(\lambda)} < 0$;  this is
easily seen to be true for sufficiently small $H$ since
$\pi_0^{(\lambda)}(0) = 1/(1-\lambda)> 0$ and
$\sum_{n=0}^{\infty}\pi_n^{(\lambda)} = 1$ for all $H$. In fact,
while the foregoing sum is absolutely convergent for $\lambda <
1/2$ it is only conditionally convergent at the Groenewold-Moyal
limit $\lambda = 1/2$, and the convergence there is very slow.)

Negative probabilities can not be dismissed as fiction. In the
present case, accepting them at face value gives the correct
expected value for the energy: differentiating the basic equation
(\ref{eq:energy}) with respect to time and setting $t = 0$ gives
$$H = \sum \pi_E\,E,$$
independent of the form of deformation. (Differentiating twice
will give the quantum second moment rather than that of the
original distribution, since the left side at $t=0$ will be
$H*H$.) In what remains we again write $\mu$ for the dimensionless
quantity $H/\hbar\omega$ when it is viewed as an ordinary scalar.

As remarked at the beginning, one view of negative probabilities
is that while they can not be observed directly we can observe
positive distributions derived from them and thereby indirectly
conclude their existence. Suppose that in our observation of the
harmonic oscillator we have a true probability distribution
$p(\mu)$ for the energy, that is, one which is non-negative,
defined for $\mu\ge 0$, and has $\int_0^{\infty}p(\mu)d\mu = 1$.
With this, the probability of the oscillator being observed at the
$n$th energy level becomes $\pi_n^{(\lambda)}(\mu,p) =
\int_0^{\infty}\pi_n^{(\lambda)}(\mu)p(\mu)d\mu.$ For some
probability distributions $p$ these will all be non-negative; such
$p$ have been called observable distributions. The existence of
sufficiently many is given immediately by (\ref{eq:fund}), for
writing
$$p_k^{(\lambda)}(\mu) =
 \frac{1}{\lambda\, k!}(\frac{\mu}{\lambda})^k
e^{-\frac{\mu}{\lambda}}$$ one has from (\ref{eq:fund}) that
$$\int_0^{\infty}\pi_n^{(\lambda)}(\mu)\,p_k^{(\lambda)}(\mu)\,d\mu =
\left\{\begin{array}{cl} \binom{k}{n}\lambda^n(1-\lambda)^{k-n} &
\mbox{ if $k\ge n$}\\0 & \mbox{0 if $k<n$} \end{array} \right. .$$
The $p_k^{(\lambda)}$ are thus observable distributions. The
coefficients on the right are the Fourier-Laguerre coefficients of
the distribution (relative to the deformation with parameter
$\lambda$). We will call the $p_k^{(\lambda)}$ \emph{basic}
observable distributions since suitable (in fact unique) linear
combinations of them (necessarily involving negative coefficients)
give all distributions with but a single non-zero Fourier-Laguerre
coefficient (i.e., we can recover those given by the individual
$\pi_n^{(\lambda)}$). The basic observable distributions are the
extreme elements of the convex cone of observable distributions
and they span that cone. Note that the basic distributions have
only a finite number of non-zero Fourier-Laguerre coefficients.

\section{The uncertainty inequality} It is a
classic computation that with the distribution $p_k^{(\lambda)}$
the mean or expected value of $\mu$, namely $\int
\mu\,p_k^{(\lambda)}(\mu)\,d\mu$, is just $\lambda(k+1)$. This is
the mean value of the energy with the basic observable
distribution $p_k^{(\lambda)}$. (Note that its minimum at $k=0$
is, as expected, just $\lambda$.) The second moment of the
distribution is $$\int \mu^2\,p_k^{(\lambda)}(\mu)\,d\mu \quad=
\quad \lambda^2(k+1)(k+2).$$ It follows that the variance is
$\lambda^2(k+1)$, the square root of which is the standard
deviation or `uncertainty'. The other evaluation of the expected
energy in $p_k^{(\lambda)}$ necessarily gives the same result, for
as observed earlier the energy calculation is independent of the
quantization (and uses the negative probabilities). In fact, we
have
$$\sum_{n=0}^\infty\,\int_0^{\infty}(n+\lambda)\pi_n^{(\lambda)}(\mu)
p_k^{(\lambda)}(\mu)\,d\mu = \sum_{n=0}^k(n+\lambda)\binom{k}{n}
\lambda^n(1-\lambda)^{k-n} = (k+1)\lambda.$$ Note that if we try
to keep the mean of the distribution constant while letting
$\lambda$ tend to zero, then the deformation becomes the normal
one and the distribution becomes a Dirac delta supported at the
mean energy. This is what was meant earlier by saying that the
normal form of quantization assumes that the energy is precisely
known.

We can now apply our selection principle.  With the basic
observable distribution $p_k^{(\lambda)}$ one has
$$\begin{aligned}\sum_{n=0}^\infty\,\int_0^{\infty}(n+\lambda)^2\pi_n^{(\lambda)}(\mu)
p_k^{(\lambda)}(\mu)\,d\mu &=
\sum_{n=0}^k(n+\lambda)^2\binom{k}{n} \lambda^n(1-\lambda)^{k-n}\\
&\qquad= (k^2+k+1)\lambda^2 +k\lambda.\end{aligned}$$ The variance
now is $k\lambda(1-\lambda)$. By our selection principle, which
asserts, in particular, that quantization should not increase
uncertainty, this quantum variance must be smaller than the
previous distribution variance. One therefore has the inequality
$$k\lambda(1-\lambda) < (k+1)\lambda^2.$$
We may not only assume that $\lambda$ is  strictly positive, but
as observed earlier, that $0 < \lambda \le 1/2$, so this implies
that
$$\lambda > \frac{k}{2k+1}.$$ This being true for all
non-negative $k$ it follows that we must have $\lambda = 1/2$,
leaving the Groenewold-Moyal form as the only one consistent with
the selection principle. Similar arguments may apply more
generally to select one quantization from a family of
cohomologically equivalent ones.

Finally, fixing $\lambda$ at $1/2$, note that the variance of the
distribution $p_k^{(1/2)}$ is $(k+1)/4$ while its quantum variance
is $k/4$, so the difference in variances is always just $1/4$. The
difference in uncertainties, however, is $(\sqrt{k+1}-\sqrt{k})/2$
which tends to 0 as the mean energy increases. This does not hold
for any $\lambda$ strictly between 0 and 1/2, so for the simple
harmonic oscillator it would be sufficient to require of a
quantization procedure that as energy tends to infinitely (and
quantization becomes unnoticeable) the difference in uncertainties
tends to zero, but this might not be a sufficiently strong
selection principle in general.

\end{document}